\documentclass[a4paper]{amsart}

\usepackage{amscd}
\newcommand{\CZ}{\mbox{$\mathbb{C}$}}
\newcommand{\QZ}{\mbox{$\mathbb{Q}$}}
\newcommand{\ZZ}{\mbox{$\mathbb{Z}$}}
\newcommand{\NZ}{\mbox{$\mathbb{N}$}}
\newcommand{\Pone}{\mbox{$\mathbb{P}_1$}}
\newcommand{\Ptwo}{\mbox{$\mathbb{P}_2$}}
\newcommand{\Pthree}{\mbox{$\mathbb{P}_3$}}

\newcommand{\name}[1]{{\sc #1}\index{#1}}

\DeclareMathOperator{\Aut}{Aut}
\DeclareMathOperator{\Pic}{Pic}

\theoremstyle{plain}
\begingroup
\newtheorem{thm}{Theorem}[section]
\newtheorem{cor}[thm]{Corollary}
\newtheorem{lem}[thm]{Lemma}
\newtheorem{prop}[thm]{Proposition}
\endgroup
\theoremstyle{remark}
\newtheorem{rem}[thm]{Remark}
\newtheorem{example}[thm]{Example}

\author{Stefan Kebekus}
\date{April 20, 1998}
\email{stefan.kebekus@uni-bayreuth.de}
\address{Stefan Kebekus\\ Mathematisches Institut der Universit\"at
Bayreuth\\ 95440 Bayreuth\\ Germany\\ FAX: +49 (0)921/55-2785 }

\title{On the Classification of 3-dimensional $SL_2(\mathbb C)$-varieties}
\thanks{The author gratefully acknowledges support by the
Graduiertenkolleg ``Komplexe Mannigfaltigkeiten'' of the Deutsche
Forschungsgemeinschaft.}

\begin{document}

\begin{abstract}
In the present work we describe 3-dimensional complex $SL_2$-varieties
where the generic $SL_2$-orbit is a surface. We apply this result to
classify the minimal 3-dimensional projective varieties with
\name{Picard}-number 1 where a semisimple group acts such that the generic
orbits are 2-dimensional.

This is an ingredient of the classification \cite{K98b} of the
3-dimensional relatively minimal quasihomogeneous varieties where the
automorphism group is not solvable.

\end{abstract}

\maketitle
\tableofcontents

\section{Introduction}

In \cite{K98b} we give a classification of the 3-dimensional
relatively minimal quasihomogeneous projective varieties where the
automorphism group is linear algebraic and not solvable. By
``relatively minimal'' we mean varieties having at most
$\QZ$-factorial terminal singularities and allowing an extremal
contraction of fiber type. These varieties always occur at the end of
the minimal model program if one starts with a projective rational
quasihomogeneous manifold whose automorphism group is not
solvable.

Certain aspects of this project utilize results on non-transitive 
$SL_2(\mathbb C)$-actions which in our opinion are of separate
interest. We have chosen to present these here as opposed to including 
them in the midst of the classification work, where the methods are
essentially different.

The aim of the first part of this paper is to describe 3-dimensional
complex $SL_2$-varieties where the generic $SL_2$-orbit is a
surface. More precisely, we give elementary criteria for the fibers of
the categorical quotient to be irreducible or normal and describe
neighborhoods of reduced fibers (see
proposition~\ref{gen:local_def}). We reduce to this case by using
concretely constructed equivariant \name{Galois} coverings which are
\'etale in codimension one. Under certain restrictions on the isotropy
group, a stronger classification is known ---see \cite{Arz98}.

In the main part of the paper we apply these results to yield the
following ingredient of the classification in \cite{K98b}.

\begin{thm}\label{S2_b21:MainThm}
Let $X$ be a $\QZ$-factorial projective 3-dimensional variety with
\name{Picard}-number $\rho (X)=1$ having at most terminal
singularities. Assume that a semisimple linear algebraic group $S$
acts algebraically on $X$ such that generic $S$-orbits are
2-dimensional. Then $X$ is isomorphic to to the smooth 3-dimensional
quadric or to one of the (weighted) projective spaces $\Pthree$,
$\mathbb P_{(1,1,1,2)}$ or $\mathbb P_{(1,1,2,3)}$.
\end{thm}

The author wishes to thank A.~T.~Huckleberry for support and many
valuable discussions. The author also would like to thank H.~Flenner
and S.~Ishii for advice on matters regarding the singularities.

\section{On the Normality of Fibers of the Categorical Quotient}

Recall that for an affine variety the quotient is defined as the
spectrum of the ring of invariant functions. The following are the
results of this section:

\begin{prop}\label{gen:fiber_is_irred}
Let $X$ be an irreducible complex affine 3-dimensional normal
$SL_2$-variety. Then all fibers of the categorical quotient map
$q:X\rightarrow Y$ are irreducible. If $X$ is additionally
\name{Cohen}-\name{Macaulay}, then a $q$-fiber is normal if it is
reduced.
\end{prop}

Under additional assumptions on the singularities, the claim is true
for non-reduced fibers as well.

\begin{prop}\label{gen:fiber_is_normal}
In the setting of proposition~\ref{gen:fiber_is_irred} assume
additionally that $X$ has at most canonical singularities. Then every
fiber of the categorical quotient is normal with it's reduced structure.
\end{prop}

Before proceeding with the proofs we recall two elementary facts:
First, the only normal affine complex $SL_2$-surfaces with
non-trivial action are
\begin{description}
\item[the smooth affine quadric $\QZ^a_2$] this space is
$SL_2$-homogeneous. The isotropy group of a point is a torus.

\item[$\Ptwo$ minus a quadric curve] this is a quotient of $\QZ^a_2$ by 
$\ZZ_2$. We denote it by $\QZ^a_2/\ZZ_2$. The isotropy group is the
normalizer of a torus.

\item[the affine cone over a rational normal curve] we denote this by
$\mathbb F^a_n$, where $n$ is the degree of the curve. The isotropy is
generated by a unipotent part and a cyclic group, isomorphic to
$\ZZ_n$. This space contains an open $SL_2$-orbit and an $SL_2$-fixed
point.
\end{description}
See \cite{Hu86} for a more detailed description. 

Second, if $X$ is a 3-dimensional $SL_2$-variety with non-trivial
action and $D_1\subset X$ is a divisor, then $SL_2$ acts non-trivially
on $D_1$. This follows directly from a linearization argument;
see~\cite[I.1.5]{H-Oe} for matters concerning linearization. In
particular, if $X$ is affine and $D_2$ is another divisor, then
$D_1\cap D_2$ must be a single point.

\begin{proof}[Proof of proposition~\ref{gen:fiber_is_irred}]
Assume without loss of generality that $\dim Y=1$, for the proposition
is trivial otherwise.  Since all $q$-fibers are connected, we must
rule out the possibility that there is a point $y\in Y$ such that
$q^{-1}(y)$ is connected and not irreducible. If this was the case,
then the irreducible components of $q^{-1}(y)$ can only meet in the unique
$SL_2$-fixed point in $q^{-1}(y)$, i.e.~$q^{-1}(y)$ is not
connected in codimension one. On the other hand, \name{Hartshorne}'s
connectedness theorem states that $X$ is connected in dimension 2 (see
\cite[Thm.~18.12 and the preceding discussion]{E95}). Now $Y$ is
normal, hence smooth, so that $q^{-1}(y)$ is \name{Cartier}. In this
situation \name{Grothendieck}'s connectedness theorem shows that
$q^{-1}(y)$ must be connected in dimension 1 (see
\cite[exp. XIII]{SGA2}), a contradiction.

If $X$ is \name{Cohen}-\name{Macaulay}, then every $q$-fiber
automatically satisfies \name{Serre}'s condition $S_2$ (see
\cite{Reid87}). If it is reduced, it's singular set is either the unique
$SL_2$-fixed point or empty. The normality follows directly from
\name{Serre}'s criterion.
\end{proof}

\begin{proof}[Proof of proposition~\ref{gen:fiber_is_normal}]
Again it is sufficient to consider the case that $\dim Y=1$. If $y\in
Y$ is a point such that $q^{-1}(y)$ has multiplicity $m>1$, let
$\Delta$ be an analytic neighborhood of $y$, isomorphic to a disk and
let $\tilde q :\tilde X \rightarrow \tilde \Delta$ be the $m$th root
fibration, associated to the restriction of $q$ to $\Delta$. If
$\tilde y$ denotes the (reduced) preimage if $y$ in $\tilde \Delta$,
then $\tilde q^{-1}(\tilde y)$ is reduced. The map $\tilde X
\rightarrow X$ is an $SL_2$-equivariant cyclic cover, branched only
over the unique $SL_2$-fixed point point in $q^{-1}(y)$, if at
all. This has two consequences: first, \cite[prop.~1.7]{Reid79}
applies, showing that $\tilde X$ has canonical singularities, so that
$\tilde X$ is \name{Cohen}-\name{Macaulay} and $\tilde q^{-1}(\tilde
y)$ is normal. Secondly, because the induced map $\tilde q^{-1}(\tilde
y)\rightarrow q^{-1}(y)$ is just the quotient by the action of the
\name{Galois} group, $q^{-1}(y)_{red}$ must also be normal.
\end{proof}

There exists a preprint of I.~V.~Arzhantsev where, using the
techniques of \cite{LV}, a proof of
proposition~\ref{gen:fiber_is_normal} is indicated for arbitrary
normal singularities.

\section{Neighborhoods of fibers}

Now we consider the neighborhood of reduced fibers.

\begin{prop}\label{gen:local_def}
In the setting of proposition~\ref{gen:fiber_is_irred}, if $Y$ is a
curve and $y\in Y$ is a point such that $q^{-1}(y)$ is reduced, then
there exists a \name{Zariski}-open neighborhood $\Delta$ of $y$ such
that $q^{-1}(\Delta)$ is equivariantly isomorphic to one of the
following:
\begin{itemize}
\item a product $\mathbb F^a_n\times \Delta$ where $SL_2$ acts on
$\mathbb F^a_n$ only
\item $\left\{ ((x,y,z),\delta)\in \CZ^3\times\Delta|
4xz-y^2=P(\delta) \right\}$, where $P\in \mathcal O(\Delta)$, having
zeros only at $y$ and $SL_2$ acts on $\CZ^3$ via the 3-dimensional
irreducible representation.
\item a quotient of the latter by $\ZZ_2$, acting with weights (1,1,1)
on $\CZ^3$ and trivially on $\Delta$.
\end{itemize}
\end{prop}

The proof follows from two technical considerations. Recall from
\cite[II.2.4]{K85} that there is an equivariant embedding
$\iota:X\rightarrow \oplus V_{k_i}$, where the $V_{k_i}$ are
irreducible $SL_2$-representation spaces.

\begin{lem}
There exists a $j\in \NZ$ such that the projection $\pi:\oplus
V_{k_i}\rightarrow V_{k_j}$ is a closed embedding if restricted to
$q^{-1}(y)$.
\end{lem}

\begin{proof}
We consider the possibilities for the central fiber separately:
\begin{description}
\item[if $\boldsymbol{q^{-1}(y)\cong \QZ^a_2/\ZZ_2}$] then every
non-trivial equivariant map is a closed embedding because the isotropy
of $\QZ^a_2/\ZZ_2$ is maximal.

\item[if $\boldsymbol{q^{-1}(y)\cong \QZ^a_2}$] the only possible images
of an $SL_2$-equivariant morphism which is not an embedding are
$\QZ^a_2/\ZZ_2$ and $\{0\}$. Both have normalizers of tori in their
isotropy groups, but $\QZ^a_2$ has not. Thus, there must be a
projection with image $\QZ^a_2$. This must be an embedding.

\item[if $\boldsymbol{q^{-1}(y)\cong \mathbb F^a_n}$] one has to rule
out that all projection map $q^{-1}(y)$ to $\{0\}$ or to $\mathbb
F^a_{kn}$, $k>1$, this being the only possible images. Assume to the
contrary and let $U<SL_2$ be a unipotent subgroup. It's fixed point
set is a line $C$, isomorphic to $\CZ$, and all projections map $C$ to
$\{0\}$ or are branched covers, ramified at zero. Thus, the rank of
the Jacobian of $\iota|_C$ drops at zero ---a contradiction to
$\iota$ being an embedding.
\end{description}
\end{proof}

Having embedded the central fibers, we show that the restriction to a
neighboring $q$-fiber is injective as well.
\begin{lem}\label{gen:local_def_lem2}
There exists a \name{Zariski}-open neighborhood $U$ of $y$ such that
for all $\eta\in U$ the restriction of $\pi$ a the $q$-fiber $X_\eta =
q^{-1}(y)$ is a closed embedding.
\end{lem}

\begin{proof}
Choose $U$ to be a maximal neighborhood of $y\in Y$ such that
$\pi(X_\eta)\not = 0$ for all $\eta\in Y$ and such that all $q$-fibers 
over $U\setminus\{y\}$ are isomorphic. Use the classification of the
2-dimensional algebraic subgroups of $SL_2$ to see that this is always 
possible. Again we perform a case-by-case check: 
\begin{description}
\item[$\boldsymbol{q^{-1}(y)\cong \mathbb F^a_1}$] in this case
$V_{k_j}$ must be $\CZ^2$. As $\pi(X_\eta)\not = \{0\}$, we have
$\pi(X_\eta)\cong \mathbb F^a_1$ and $X_\eta$ must be isomorphic to
$\mathbb F^a_1$ itself, there being no $SL_2$-equivariant cover.

\item[$\boldsymbol{q^{-1}(y)\cong \mathbb F^a_2}$] here $V_{k_j}$ is the 
irreducible 3-dimensional representation space. The only $SL_2$-invariant
divisors in here are $\mathbb F^a_2$ and smooth quadrics. Arguing as
above, one must show that the generic $q$-fiber $X_\eta$ is not
isomorphic to a cover of $\mathbb F^a_2$ or $\QZ^a_2$, i.e. $X_\eta
\not\cong \mathbb F^a_{1}$. If this was the case, then linearize the
center $Z$ of $SL_{2}$ at a smooth point of $q^{-1}(y)$. This gives an
analytic curve germ $C\subset X$, invariant under $Z$ and intersecting
$q^{-1}(y)$ transversally in a single point. As $Z$ is not contained
in the isotropy group of any point in $X_\eta$ other than 0, $C$ must
intersect the neighboring fiber twice.  This is a contradiction to
$q^{-1}(y)$ being reduced.

\item[$\boldsymbol{q^{-1}(y)\cong \mathbb F^a_n}$ where
$\boldsymbol{n=3}$ or $\boldsymbol{n>4}$]a similar linearization
argument as above, using a $\ZZ _{n}$ from the isotropy group of a
generic point in $q^{-1}(y)$, shows that the generic $X_\eta$ must
contain a $\ZZ _{n}$-fixed curve. Classification yields that $X_\eta
\cong \mathbb F^a_{kn}$ for one $k\in \NZ $. But $k$ must be 1: every 
$X_\eta$ contains a curve which is $\ZZ _{kn}$-fixed and $q^{-1}(y)$
must, too.

\item[$\boldsymbol{q^{-1}(y)\cong \mathbb F^a_4}$] here $V_{k_j}$ is the
irreducible 5-dimensional representation space where the only
$SL_2$-invariant surfaces are $\mathbb F^a_4$ or are isomorphic to
$\QZ^a_2/\ZZ_2$. The linearization argument used above rules out that
$X_\eta\cong \QZ^a_2$ or a cover of $\mathbb F^a_4$.

\item[$\boldsymbol{q^{-1}(y)\cong \mathbb Q^a_2}$]here $V_{k_{j}}$ 
contains two types of 2-dimensional $SL_{2}$-orbits: $\mathbb Q^a_2$
and $\mathbb F^a_{k_{j}}$. We know that $\pi (X_\eta)\cong \mathbb
Q^a_2$, as otherwise $\pi(X_\eta)$ must contain a $U$-pointwise fixed
curve and $\pi(q^{-1}(y))$ must, too. A contradiction. Again $\pi
|_{X_\eta}$ must be injective as there is no $SL_{2}$-equivariant
cover of $\mathbb Q^a_2$.

\item[$\boldsymbol{q^{-1}(y)\cong \mathbb Q^a_2/\ZZ _{2}}$] apply 
the linearization argument involving a generic isotropy group,
i.e. the normalizer of a torus to see that the neighboring $q$-fibers
cannot be isomorphic to $\mathbb Q_{2}$. Now argue as in the last
case.
\end{description}
\end{proof}

With this information we start the 

\begin{proof}[Proof of proposition~\ref{gen:local_def}]
Choose $\Delta\subset Y$ as in lemma~\ref{gen:local_def_lem2}. Then the
map $(\pi\circ \iota)\times q:X\rightarrow V_{k_j}\times Y$ is
injective if restricted to $q^{-1}(\Delta)$.

Recall that every irreducible representation space of $SL_2$ contains
a unique $SL_2$-orbit whose closure is isomorphic to $\mathbb 
F^a_n$. Thus the claim of proposition~\ref{gen:local_def} holds if all
$q$-fibers are isomorphic to $\mathbb F^a_n$.

If $q^{-1}(y)\cong \mathbb F^a_2$, and the generic fiber is a smooth
quadric, then $q^{-1}(\Delta)$ can be equivariantly embedded into
$V_2\times \Delta$. Equip $V_2$ with coordinates $(x,y,z)$, fix a
torus $T<SL_2$ and note that it's fixed point set $V_2^T$ is a
1-dimensional linear subspace. If $y$ is the linear coordinate on
$V_2^T$, then the intersection $X^T:= X\cap (V_2^T\times\Delta)$ is
given by $\{-y^2=P(\delta)\}$ where $P\in\mathcal O(\Delta)$. This is
because $X^T$ is 2:1 over $\Delta$ and invariant under multiplication of
$V_2$ with -1. By choice of $\Delta$, $P$ has no zero on
$\Delta\setminus \{0\}$. Now $X$ being uniquely determined by $X^T$ as
$X=\overline{SL_2.X^T}$ shows that $X$ is given by
$\left\{((x,y,z),\delta)\in \CZ^3\times\Delta| 4xz-y^2=P(\delta)
\right\}$, all $SL_2$-invariant surfaces in $V_2$ being given as
$4xz-y^2=const$ after proper choice of coordinates. Thus, the claim is
shown as well.

If all fibers are isomorphic to $\QZ^a_2$ and $V_{k_j}\not\cong V_2$,
then argue similarly: $V_{k_j}^T$ is 1-dimensional and
$X^T=X\cap(V_{k_j}^T\times\Delta)$ is given by $\{-y^2=P(\delta)\}$ where
$P\in \mathcal O^*(\Delta)$. We show that $X$ is isomorphic to
$\left\{((x,y,z),\delta)\in \CZ^3\times\Delta| 4xz-y^2=P(\delta)
\right\}=:X_2\subset V_2\times\Delta$. A linear identification of
$V_2^T$ and $V_{k_j}^T$ yields an isomorphism between $X^T$ and
$X_2^T$. Let $\Gamma^T \subset (V_2^T\times\Delta) \times
(V_{k_j}^T\times\Delta) \subset (V_2\times\Delta) \times
(V_{k_j}\times\Delta)$ be the graph and set
$\Gamma:=SL_2.\Gamma^T\subset (V_2\times\Delta) \times
(V_{k_j}\times\Delta)$. Now $X^T$ and $X_2^T$ both having isotropy
group $T$ at any point implies that $\Gamma$ is the graph of a
bijective morphism, i.e.~an isomorphism between the (normal) varieties
$X$ and $X_2$.

If $q^{-1}(y)\cong \mathbb F^a_4$ or all $q$-fibers are isomorphic to
$\QZ^a_2/\ZZ_2$ one uses the same line of argumentation with the only
difference that the $SL_2$-invariant surfaces in the 5-dimensional
representation space are given by the ideal 
\begin{align*}
3d^2-8ce+4\delta e, && cd-6be+\delta d \\
3bd-48ae+2\delta c+2\delta^2, && c^2-36ae+2\delta c+\delta^2 \\
bc-6ad+\delta b, && 3b^2-8ac+4\delta a.
\end{align*}
This variety is a quotient of $\{4xz-y^2=\delta\}$ by $\ZZ_2$, where
the $SL_2$-equivariant quotient map is given by $(x,y,z)\mapsto
(x^2,2xy,2xz+y^2,2yz,z^2)$.
\end{proof}

The next lemma covers a special case which we will need to consider
later.

\begin{prop}\label{gen:glob_def}
In the setting of proposition~\ref{gen:local_def}, if $Y\cong \CZ$ and
all fibers over $Y\setminus\{0\}$ are isomorphic, then $X$ is
equivariantly isomorphic to
\begin{itemize}
\item $\mathbb F^a_n\times \CZ$ where $SL_2$ acts on
$\mathbb F^a_n$ only
\item $X_k:=\left\{((x,y,z),\delta)\in \CZ^3\times\CZ| 4xz-y^2=\delta^k
\right\}$ where $k\in \NZ$ and $SL_2$ acts on $\CZ^3$ via the 3-dimensional
irreducible representation.
\item a quotient of the latter by $\ZZ_2$, acting with weights (1,1,1)
on $\CZ^3$ and trivially on the base.
\end{itemize}
\end{prop}
\begin{proof}
If all $q$-fibers are isomorphic to $\mathbb F^a_n$, then
proposition~\ref{gen:local_def} shows the local triviality. Note that
the only automorphisms of $\mathbb F^a_n$ commuting with the
$SL_2$-action are in $\CZ^*$. But $H^1(\CZ,\mathcal O^*)$ is trivial
so that the local trivializations glue together to give a global one.

If the generic fiber is $\QZ^a_2$, then employ the same methods as in
the proof of proposition~\ref{gen:local_def}: embed $X\rightarrow \oplus
V_{k_i}\times \CZ$ and assume that $(\pi_0\times Id):\oplus
V_{k_i}\times\CZ\rightarrow V_{k_0}\times\CZ$ is an embedding, if
restricted to $q^{-1}$ of a neighborhood of $y$. Choose a torus $T$
and let $X^T\subset V_{k_0}^T$ be the $T$-fixed point set. Assuming
without loss of generality that $y=0$, $(\pi_0\times Id)(X^T)$ is
given as $\{-y^2=c\cdot\delta^k\cdot\prod(y_j-\delta)^{m_j}\}$ where
$y_j\not = 0$, $\delta$ is to coordinate on $Y\cong \CZ$ and $c\not =
0$ is a constant. We know that $X^T$ is locally (analytically)
reducible over each of the $y_j$. It's $(\pi_0\times Id)$-image is,
as well. Thus, the $m_j$ are even.

Let $U_0\subset \CZ$ be the maximal set such that all fibers are
isomorphic to $\QZ^a_2$. To construct an isomorphism $X_k\rightarrow
X$ over $U_0$, it is necessary to find an isomorphism between $X^T\cap
q^{-1}(U_0)$ and $X_k^T:=\{-y^2=\delta^k|\delta\in U_0\}$ and then
apply the construction from the proof of
proposition~\ref{gen:local_def}, involving the graph $\Gamma$. Note
that $X^T\cap q^{-1}(U_0)$ and $X_k^T$ are both smooth and have a
birational morphism onto $(\pi_0\times Id)(X^T)$, the latter being
given by
\begin{align*}
X_k^T &\rightarrow (\pi_0\times Id)(X^T) \\
(y,\delta) & \mapsto \left(y\prod(y_j-\delta)^{\frac{m_j}{2}}\sqrt{c}),\delta\right).
\end{align*}
Thus, they must be isomorphic. Now the construction gives an
isomorphism over $U_0$.

If $U_0\subset \CZ$ is not the whole of $\CZ$, then set
$U_1:=\CZ\setminus\{y_1,\ldots,y_k\}$. Recall that $V_{k_0}$ is
necessary 3-dimensional, equip it with coordinates $x$, $y$ and $z$
and set
\begin{align*}
\left\{4xz-y^2=\delta^k \right\} & \rightarrow \left\{4xz-y^2=\delta^k\cdot\prod(y_j-\lambda)^{m_j}\right\} \\
((x,y,z),\delta) & \mapsto \left(\prod
(y_j-\delta)^{\frac{m_j}{2}}(x,y,z),\delta\right)
\end{align*}
where $((x,y,z),\delta)$ are coordinates on $\CZ^3\times\CZ$.

We have to show that the two local isomorphisms over $U_0$ and $U_1$
agree. Note that the only automorphisms of $\QZ^a_2$ commuting with the
$SL_2$-action are in $\ZZ_2$. Now $H^1(\CZ,\ZZ_2)$ being trivial shows
hat after multiplying one of the local isomorphisms with (-1), if
necessary, we can always glue.

Again the analogous construction works if the generic fiber is
isomorphic to $\QZ^a_2/\ZZ_2$.
\end{proof}

\begin{rem}
Propositions~\ref{gen:local_def} and \ref{gen:glob_def} could also be
proved using elementary deformation theory, see~\cite{Pi74}.
\end{rem}

Now we describe certain quasi-projective varieties which will play an
important role in the next chapter. For this, a categorical quotient
of a quasi-projective variety is an invariant affine surjective
morphism $q:X\to Y$ which is a categorical quotient on an affine cover 
of the base.

\begin{prop}\label{gen:toric}
Let $X$ be a quasi-projective normal complex $SL_2$-variety with at
most terminal singularities and categorical quotient $q:X\to
\Pone$. If every $q$-fiber over $\CZ^*\subset \Pone$ is isomorphic to
$\CZ^2$, then the singularities of $X$ are of type
$\frac{1}{n}(1,1,-1)$ and $\frac{1}{m}(1,1,-1)$ (i.e.~are locally
isomorphic to $\CZ^3/\ZZ_n$, where $\ZZ_n$ acts with weights (1,1,-1))
and $X$ is toric. Here $n$ and $m$ are the multiplicities of the
exceptional $q$-fibers.
\end{prop}
\begin{proof}
Set $X^*:=q^{-1}(\CZ^*)$. By proposition~\ref{gen:local_def}, $X^*$ is
a locally trivial $\CZ^2$-bundle. Since the transition functions must
commute with $SL_2$, they are in $\mathcal O^*$, and $X^*$ must be the
sum of two line bundles. However, $\Pic(\CZ^*)=0$, so that $X^*$ is
isomorphic to the trivial bundle. Choose two different unipotent
subgroups $U_1$, $U_2<SL_2$ and two sections $\sigma_1$ and $\sigma_2$
in $X^*$ over $\CZ^*$ which do not have zeros and such that $\sigma_i$
is pointwise $U_i$-fixed. The $\sigma_i$ yield an isomorphism
$X^*\rightarrow \CZ^2\times\CZ^*$. Let $(\CZ^*)^3$ act on $X^*$ in
these coordinates by 
$$ 
(r,s,t)((x,y),\lambda)\mapsto ((rx,sy),t\lambda) 
$$ 

Set $X^0:=q^{-1}(\CZ)$ and let $n$ be the multiplicity of
$q^{-1}(0)$. Let $\gamma:\tilde X^0\rightarrow X^0$ be the $n$th root
fibration. By proposition~\ref{gen:glob_def}, $\tilde X^0$ is the
trivial $\CZ^2$-bundle over $\CZ$. Choosing sections $\tau_i$ as above
yields coordinates $((\tilde x,\tilde y),\lambda)$ on $\tilde
X^0$. Then the closures of the pull-back of the $\sigma_i$ are given
by $\overline{\gamma^{-1}(\sigma_1)}=\{x=0, y=\lambda^{k_1}\}$ and
$\overline{\gamma^{-1}(\sigma_2)}=\{x=\lambda^{k_2}, y=0\}$ and the
$(\CZ^*)^3$-action on $X^*$ pulls back to $(r,s,t)((\tilde x,\tilde
y),\lambda)\mapsto ((t^{k_2}rx,t^{k_1}sy,t\lambda)$.  In particular,
the $(\CZ^*)^3$-action can be extended to the whole of $\tilde X^0$.

By construction, $X^0$ is a cyclic quotient of $\tilde X^0$ by the
group $\ZZ_n$. If $(f,\lambda)\in \CZ^2 \times \CZ$ and $\xi\in
\ZZ_{n}$ is a primitive root, then we may write without loss of
generality $\xi(f,\lambda)=(a(\xi,\lambda)f,\xi\cdot \lambda)$,
where $a(\xi,\lambda)\in \Aut (\mathbb F^a_n)$, commuting with $SL_2$,
i.e.~$a(\xi,\lambda)\in\CZ^*$. Since there is no non-constant
algebraic morphism from $\CZ$ to $\CZ^*$, $a(\xi,\lambda)$ does not
depend on $\lambda$. We may thus view $\ZZ_n$ as acting on $\CZ^2 \times
\CZ$ with weights $(a,a,1)$.  Note that the $\ZZ_n$-quotient
map commutes with the action of $(\CZ^*)^3$, i.e.~$X^0$ is
toric. A similar argument over $\Pone\setminus \{0\}$ shows that $X$
is toric.

In order to show that $X$ has singularities of type
$\frac{1}{n}(1,1,-1)=\frac{1}{n}(-1,-1,1)$, i.e.~that $a=1$, note that
the $a$ and $n$ must be coprime (this is because $\gamma$ is \'etale
in codimension one). The classification theory of terminal
singularities asserts that only $a=-1$ is possible if the quotient is
supposed to be terminal, see~\cite[sect.~5.1]{Reid87}. Again, the same
argument applies to $q^{-1}(\Pone\setminus \{0\})$.
\end{proof}

\section{An Application}

Our main goal here is to prove theorem~\ref{S2_b21:MainThm}, but first
we construct the example which is actually realized.

\begin{example}\label{S2_b21:Type_I_example}
Let $X$ be the weighted projective space $\mathbb P_{(1,1,2,3)}$. By
\cite[p. 35]{Ful93}, identify $X$ with the  toric variety $X$ whose
fan is constructed from the vectors $e_{1}$, $e_{2}$, $e_{3}$ and
$v=(-1,-2,-3)\in \ZZ ^{3}$, i.e.~whose cones are generated by any
three of these vectors.  We claim that $X$ has $\QZ $-factorial
terminal singularities and \name{Picard-}number $\rho
(X)=1$. Furthermore, $SL_2$ acts on $X$ such that the generic orbit is
2-dimensional.

\begin{description}
\item[$\boldsymbol X$ has $\QZ$-factorial terminal singularities] The 
cones generated by $(e_{1,}e_{2},e_{3})$ and $(e_{2},e_{3},v)$
describe smooth varieties. The cone generated by $(e_{1},e_{2},v)$ can
be brought into the form $(e_{1},e_{2},(-1,-2,3))$ using the matrix
$Diag(1,1,-1)\in GL(2,\ZZ )$. The latter cone is known
(see \cite[p. 35]{Ful93}) to describe the singularity $\CZ
^{3}/\ZZ_3$, where $\ZZ _{3}$ acts with weights
$(1,1,-1)$. Analogously, 
$$
\begin{pmatrix}
1&0&0\\0&1&0\\0&-1&1
\end{pmatrix}
(e_{1},e_{3},v) =(e_{1},e_{3},(-1,2,-1)).
$$
This describes the singularity $\CZ^3/\ZZ _{2}$, where $\ZZ _{2}$
acts with weights $(1,1,1)$. Both singularities are terminal
---compare \cite[thm. on p. 379]{Reid87}.

\item[$\boldsymbol X$ has \name{Picard}-number one] By \cite[p. 64f]{Ful93}, 
$\Pic(X)=\ZZ $.

\item[$\boldsymbol{SL_2}$-action] If $[x:y:z:w]$ are weighted
homogeneous coordinates associated to the weights $(1,1,2,3)$, let
$SL_2$ act on $x$ and $y$ via the 2-dimensional representation.
\end{description}
\end{example}

Since the proof of theorem~\ref{S2_b21:MainThm} is rather long,
we subdivide it into a number of lemmata. For the rest of this paper, 
we use the notation of theorem~\ref{S2_b21:MainThm} without
mentioning it further.

\begin{lem}
In the situation of theorem~\ref{S2_b21:MainThm}, $S\cong SL_2$.
\end{lem}

\begin{proof}
If it was not, then the generic $S$-orbit is a homogeneous divisor.
Two of them don't intersect, a contradiction to $\rho (X)=1$.
\end{proof}

\begin{lem}\label{S2_b21:S-stab-div}
Every irreducible $S$-invariant divisor $D_{\alpha }\subset X$ is
$S$-quasihomogeneous and it's normalization is isomorphic to either
\begin{enumerate}
\item the projective cone over the rational normal curve of degree
$n$, $\mathbb F_{n}$, where $S$ acts with a fixed point (this includes
$\Ptwo$),
\item $\Ptwo$, where $S$ acts via the 3-dimensional irreducible
$SL_2$-representation, or
\item the \name{Hirzebruch} surface $\Sigma _{0}$, where $S$ acts diagonally
\end{enumerate}
There exists a curve $C\subset X$, $C\cong\Pone$ such that
(set-theoretically) $C=D_\alpha\cap D_\beta$ for all $S$-invariant
divisors $D_\alpha$ and $D_\beta$. Furthermore, every irreducible
$S$-invariant divisor is locally (analytically) irreducible at any point
of $C$.
\end{lem}

\begin{proof}
We have remarked at the very beginning that $D_{\alpha }$ cannot be
pointwise $S$-fixed. Suppose that the generic $S$-orbit in $D_\alpha$
was 1-dimensional. Then $D_\beta$ intersects $D_\alpha$ in finitely
many orbits and therefore has empty intersection with a generic
$S$-orbit in $D_{\alpha}$, a contradiction to $\rho
(X)=1$. Consequently, $D_{\alpha }$ is $S$-quasihomogeneous and if
$\eta :\tilde D _{\alpha }\rightarrow D_{\alpha }$ is the
normalization, then classification (see \cite{Hu86}) shows that
$\tilde D _{\alpha }$ must be isomorphic a variety in the list, or to
$\Sigma _{n}$, $n>0$ with $S$ stabilizing two sections.

Due to $\rho (X)=1$ there is a number $k\in \NZ$ such that $kD_{\beta
}$ is \scshape Cartier \upshape and ample. In particular, $\eta
^{*}(kD_{\beta })$ is an effective ample divisor with $S$-invariant
support. This is possible in all cases save $\Sigma_n$, $n>0$. In the
allowed cases, there is a unique irreducible $S$-invariant curve
$C\subset \tilde D _{\alpha }$ which yields the assertion.
\end{proof}

One of the key points in the proof of theorem~\ref{S2_b21:MainThm} is
the following local description of the $D_\alpha$ in the neighborhood
of $C$.

\begin{lem}\label{S2_b21:local_situa}
Let $ x\in C$ and $ T<S$ be a torus fixing $ x$. Then there exists a
2-dimensional $ T$-representation space $E$ with positive weights, a
neighborhood $ V$ of $ 0\in E$, a neighborhood $ \Delta $ of $ 0\in
\mathbb C$ and an immersion $ \phi :V\times \Delta \rightarrow
U\subset X$ with the following properties:
\begin{enumerate}
\item $ \phi ^{-1}(C)=\{0\}\times \Delta $
\item $ D_{\alpha }\subset X$ is an $ S$-invariant divisor if and only if
there exists a $ T$-invariant curve $ N\subset E$ such that $ \phi
^{-1}(D_{\alpha })=N\times \Delta $. 
\end{enumerate}
\end{lem}

\begin{proof}
Let $ U\subset T_{x}X$ be a sufficiently small neighborhood of $ 0$
and let $ \lambda : U\rightarrow X$ be a linearization of the $
T$-action. As $ C$ is smooth, the tangent space $ T_{x}C\subset
T_{x}X$ coincides with one of the weight spaces. Take two weight
vectors $ v_{1}$ and $ v_{2}\in T_{x}X$ which, together with $ T_{x}C$
span the whole space $ T_{x}X$. Let $ E$ be the space spanned by $
v_{1}$ and $ v_{2}$.

As a first step, we claim that after replacing $ T$ be $ T^{-1}$, if
necessary, all weights of the $ T$-action on $ T_{x}X$ are
positive. In order to show this it is sufficient to show that for
generic $ y\in U$ the $ T$-orbit $ Ty$ contains $ x$ in the
closure. This, however, is true because $ \lambda (y)$ is contained in
an $ S$-invariant divisor and $x\in \overline{Ty}$ by the classification
of lemma \ref{S2_b21:S-stab-div}.

In order to construct the map $ \phi $, note that $ \lambda |_{V}$ is
immersive. The image of the tangential map $ T(\lambda |_{V})$ is
transversal to $ T_{x}C$. Let $ H<S$ be a unipotent one-parameter
group not fixing $ x$. The associated vector field, evaluated at $ x$
is contained in $ T_{x}C$ so that the map
$$
\begin{array}{cccc}
\phi : & V\times H & \rightarrow  & X\\
 & (v,h) & \mapsto  & h\cdot \lambda (v)
\end{array}
$$
has maximal rank at $ (0,0)$. Thus, $ \phi $ is invertible in a small
neighborhood. Property (1) holds by construction. 

In order to show property (2) it is sufficient to consider irreducible
$D_{\alpha }$. Claim that $D_{\alpha }$ is $S$-invariant if and only if
there exists a point $ y\in D_{\alpha }\cap U$ such that $D_{\alpha
}=\overline{H.\overline{T.y}}$. Indeed, if $ D_{\alpha }$ is
$S$-invariant, then it contains $ C$, and therefore also a point $ y\in
D_{\alpha }\cap \lambda (V)\setminus C$ and $\overline{T.y}$ is
necessarily a curve containing $ x$ in the closure. Therefore
$\overline{T.y}$ is \em not \em $ H$-invariant and
$\overline{H.\overline{T.y}}$ is an irreducible component of
$D_{\alpha }$, hence equal to $ D_{\alpha }$. This shows already that
if we set $N:=\overline{T.\phi ^{-1}(y)}$, then $\phi (N\times
\Delta )$ is contained in $D_{\alpha }$. If $\phi ^{-1}(D_{\alpha
})\not =N\times \Delta$, then it contains another irreducible
component, a contradiction to the local irreducibility of $ D_{\alpha
}$.  
\end{proof}

We utilize the local description to draw conclusions concerning the
global configuration of the divisors $D_\alpha$.

\begin{cor}\label{S2_b21:coprimeness}
There are at least two different $S$-invariant divisors $D_{0}$ and
$D_{\infty }$ in $ X$ which are smooth along $C$. Unless $X$ is
isomorphic to the smooth 3-dimensional quadric $\QZ_3$, to $\Pthree$
or to $\mathbb P_{(1,1,1,2)}$, the normalizations are isomorphic to
$\tilde D _{0}\cong \mathbb F_{n}$ and $\tilde D _{\infty }\cong
\mathbb F_{m}$ with $n,m>1$. If $\tilde D_\alpha$ is the normalization
of a generic $S$-invariant divisor, then either
\begin{enumerate}
\item $m$ and $n$ are coprime and $\tilde D_\alpha \cong \Ptwo$ where $S$ 
acts with a fixed point, or

\item $m$ and $n$ are even, $\frac{m}{2}$, $ \frac{n}{2}$ are coprime and 
$\tilde D_\alpha \cong \Sigma_0$, or

\item $m$ and $n$ are divisible by four, $\frac{m}{4}$, $\frac{n}{4}$ are 
coprime and $\tilde D_\alpha \cong \Ptwo$ where $S$ acts via the
3-dimensional irreducible representation.
\end{enumerate}
\end{cor}

\begin{proof}
Taking $ N$ to be one of the weight spaces in $E$, lemma
\ref{S2_b21:local_situa} immediately yields $D_{0}$ and $D_{\infty
}$. Note that by lemma \ref{S2_b21:S-stab-div} two $S$-invariant
divisors intersect in $C$ only so that the generic $S$-invariant
divisor $D_{\alpha }$ does not meet the singular set of $X$. Use the
standard argument linearizing the $S$-action at a fixed point to
exclude the possibility that $D_{\alpha }\cong \mathbb F_{n}$ where
$n>1$. Thus, $ D_{\alpha }$ is smooth away from $C$.

Secondly, remark that if $X$ is a cone then the classification from
\cite[thm.~3.3 and cor.~3.4]{Mori82} yields that $X\cong\Pthree$ or
$\mathbb P_{(1,1,1,2)}$ if $X$ is assumed to have $\QZ$-factorial and
terminal singularities; note that a cone over $\Sigma_0$ is never
$\QZ$-factorial as there are \name{Weyl}-divisors intersecting in
a single point.

Recall that $X$ is isomorphic to a cone or to $\QZ_3$ if there is a
\name{Cartier} divisor in $X$ which is isomorphic to $\mathbb P_{2}$
or $\Sigma _0$; see \cite[thms.~1 and 5]{Bad82} for the cases that $X$
is smooth or that $D_\alpha\cong
\Ptwo$ and \cite[thm.~3]{Bad84} for the remaining case.

Thus, excluding this case amounts to saying that the normalizations of
$D_0$ and $D_\infty$ are isomorphic to $\mathbb F_{\bullet}$, since
otherwise $D_\bullet \setminus C$ would be homogeneous, would not
intersect the (finite) singular set of $X$ and would thus be
\name{Cartier}. The indices $n$ and $m$ are exactly the weights of the
$ T$-action on $ E$, as given by lemma~\ref{S2_b21:local_situa}. The
possible weights of the $T$-action on the $SL_2$-quasihomogeneous
surfaces $D_{\alpha }$ (see the classification of
lemma~\ref{S2_b21:S-stab-div}) and the fact that $N$ {\em must} be
singular give conditions (1)--(3).
\end{proof}

Note that the set of semi-stable points with respect to the unique
lifting of the $SL_2$-action to $\mathcal{O}(D_\alpha)$ is $X\setminus
C$. Let $q:X\setminus C \to Y$ denote the resulting quotient in the
sense of geometric invariant theory.

\begin{cor}\label{S2_b21:multiplicities}
We have $Y\cong \Pone$ and either $X\cong \QZ_3$, $\Pthree$ or
$\mathbb P_{(1,1,1,2)}$ or and there are points $0$, $\infty\in \Pone$
such that $q^{-1}(0)=n'D_0$, $q^{-1}(\infty)=m'D_\infty$ and all other
$q$-fibers are reduced. Here $n'=n$, $\frac{n}{2}$ or $\frac{n}{4}$,
according to the cases of corollary~\ref{S2_b21:coprimeness}; $m'$
similarly.
\end{cor}
\begin{proof}
The description of lemma~\ref{S2_b21:local_situa} guarantees that the
quotient map extend to a rational map $X\dasharrow Y$ which becomes
regular if we perform a weighted blow-up of $C$ with weights $n$ and
$m$. Since the exceptional set of this blow-up is rational, $Y$ is, as 
well. Thus $Y\cong \Pone$. In particular, all $q$-fibers are linearly
equivalent, and the $D_\alpha$ are linearly equivalent up to positive
multiplicities. 

In order to see that all the $D_\alpha$ have multiplicity 1 as
$q$-fibers, it is sufficient to see that the divisors $D_\alpha$ are
linearly equivalent. By lemma~\ref{S2_b21:local_situa}, $D_\alpha$ is
locally given by a curve $N$ having the equation $x^{(n')}=y^{(m')}$,
$D_{0}=\{x=0\}$ and $D_{\infty }=\{y=0\}$. Thus,
$$ 
D_{0}.D_{\alpha } = m'C, \quad D_{1}.D_{\alpha } = n'C, \quad
D_{\alpha }.D_{\alpha } = n'm'C.  
$$ 
Consequently $$ D_{0} \sim
\frac{1}{n'}D_{\alpha } \quad\text{and}\quad D_{\infty } \sim
\frac{1}{m'}D_{\alpha }  
$$ 
as $\QZ$-divisors. This finishes the proof.
\end{proof}

\subsection{Proof the Theorem~\ref{S2_b21:MainThm}}

With these preparations we start the proof the main
theorem~\ref{S2_b21:MainThm}. If $X\cong \mathbb{Q}_3$, $\Pthree$ or
$\mathbb P_{(1,1,1,2)}$, we can stop here. Otherwise, we are in one of
the cases (1)--(3) or corollary~\ref{S2_b21:coprimeness}. We treat
these cases separately.

\begin{proof}[Proof of~\ref{S2_b21:MainThm} in case (1) of
corollary~\ref{S2_b21:coprimeness}] 
By proposition~\ref{gen:local_def}, all $q$-fibers over
$\Pone\setminus\{0,\infty\}$ are isomorphic to $\CZ^2$. By
proposition~\ref{gen:toric}, $X$ has two singularities of type
$\frac{1}{n}(1,1,-1)$ and $\frac{1}{m}(1,1,-1)$, and $X\setminus C$ is
toric. Since $X$ is smooth along $C$, the associated vector fields
extend to $X$, showing that $X$ is toric, too.

Consequence: $X$ can be given as a fan in $\ZZ^3$. Let $\sigma _{1}$
and $\sigma _{2}\subset \ZZ^3$ be the cones describing the smooth
$(\CZ^*)^3$-fixed points on $C$, $\sigma _{3}$ describe the point
$\frac{1}{n}(1,1,-1)$ and $\sigma _{4}$ be associated to
$\frac{1}{m}(1,1,-1)$. There can be no further fixed points. 

Choose coordinates such that $\sigma _{1}$ is spanned by the unit
vectors $(e_{1},e_{2},e_{3})\in \ZZ ^{3}$. Because every cone is
spanned by four rays, there must be a vector $v=(a,b,c)\in \ZZ ^{3}$
such that $\sigma _{2}$, $\sigma _{3}$ and $\sigma _{4}$ are spanned
by 2 unit vectors and $v$ each. After renaming the $e_{i}$, if
necessary, assume that $\sigma _{2}=(e_{1},e_{2},v)$, $\sigma
_{3}=(e_{1},e_{3},v)$ and $\sigma _{4}=(e_{2},e_{3},v)$. We will find
out the possibilities for $v$.  First, note that two cones must not
intersect in anything but a face.  Thus, $a$, $b$ and $c$ must be
negative. We use the local description of the singularities:

\begin{description}
\item[$\boldsymbol{\sigma_2}$ is smooth] consequently, $(e_{1},e_{2},v)$ 
must be a basis of $\ZZ ^{3}$ and $c=-1$.

\item[$\boldsymbol{\sigma_3}$ is $\boldsymbol{\frac{1}{n}(1,1,-1)}$] It 
is known (see \cite[p. 35]{Ful93}) that the cone generated by
$(e_{1},e_{3},-(n-1)e_{1}+ne_{2}-e_{3})$ corresponds to a singularity
of type $\frac{1}{n}(1,1,-1)$. Thus, there exists a $g\in GL(3,\ZZ )$
such that
$g(e_{1},e_{3},-(n-1)e_{1}+ne_{2}-e_{3})=(e_{1},e_{3},(a,b,-1))$. Calculating
the product $$
\begin{pmatrix}
1&\alpha&0\\
0&\beta&0\\
0&\gamma&1
\end{pmatrix}
\begin{pmatrix}
-(n-1)\\
n\\
-1
\end{pmatrix}
=
\begin{pmatrix}
-n+1+\alpha n\\
\beta n\\
\gamma n -1
\end{pmatrix}
\begin{matrix}
!\\
=\\
\ \\
\end{matrix}
\begin{pmatrix}
a\\
b\\
-1
\end{pmatrix}
$$
yields $\gamma =0$, and $a\in \ZZ n+1$. Since $\det g\in \pm
1$, $\beta \in \pm 1$. The inequality $b<0$ gives $b=-n$.

\item[$\boldsymbol{\sigma_4}$ is $\boldsymbol{\frac{1}{m}(1,1,-1)}$] 
Similar to the above there is a $g\in GL(3,\ZZ )$ such that
$g(e_{2},e_{3},me_{1}-(m-1)e_{2}-e_{3})=(e_{1},e_{3},(a,-n,-1))$. The
same calculation shows $a=-m$ and $b\in \ZZ m+1$.
\end{description}

Summarizing the above, we need to find all $n$ and $m$ such that there
are numbers $\mu,\nu\in \ZZ $ with
\begin{align}
\label{S2_b21:eq1} m &=\mu n-1 \\
\label{S2_b21:eq2} n &=\nu m-1
\end{align}

By assumption, $n$ and $m$ are coprime so that we can always assume
without loss of generality that $n>m>1$.  Then
equation~\ref{S2_b21:eq1} holds iff $\mu=1$ and $m=n-1$. Inserting
this into equation~\ref{S2_b21:eq2} gives $m(\nu-1)=2$ which in turn
implies $m=2$. Now compare $v$ to the description of
example~\ref{S2_b21:Type_I_example}.
\end{proof}

\begin{proof}[Proof of~\ref{S2_b21:MainThm} in case (2) of
corollary~\ref{S2_b21:coprimeness}] Let us begin by giving a detailed
description of this case over a trivialization.  Set $X^0:=X\setminus
D_\infty$. By corollary~\ref{S2_b21:multiplicities}, $q^{-1}(0)$ has
support on $D_0$ and multiplicity $n'$; this is the only $q$-fiber
with non-trivial multiplicity over $\CZ$. Let $\tilde q:\tilde
X^0\rightarrow \CZ$ be the $n'$th root fibration associated to
$q:X^0\rightarrow \CZ$. Now $\tilde X^0$ is a quotient of $X_k$ by the
cyclic group $\ZZ_{n'}$ acting freely in codimension 1. Choose an
analytic disk $\Delta\subset\CZ$ around 0 which is
$\ZZ_{n'}$-invariant. Then, after proper choice of coordinates,
$\tilde q^{-1}(\Delta)\cong\{4xz-y^2=\delta^k\}$ as ensured by
proposition~\ref{gen:local_def}.

It is elementary to see that every automorphism of $\tilde
q^{-1}(\Delta)$ over $\Delta$ commuting with $SL_2$ is given by
$((x,y,z),\delta)\rightarrow (\pm t^{\frac{k}{2}}(x,y,z),t\delta)$ for
some $t\in\CZ^*$. Thus, the action of $\ZZ_{n'}$ on $X_k$ extends to
$\CZ^3\times\Delta$ where $\tilde q^{-1}(\Delta)$ is
$SL_2$-equivariantly embedded. Since the action of $\ZZ_{n'}$ must
commute with $SL_2$, the weights of the $\ZZ_{n'}$-action on
$\CZ^3\times\{0\}$ must be equal.

An analogous construction can be given at $\infty$. We now show that
$n'$ and $m'$ are not coprime. This contradicts the assumption.

Consider the following cases:
\begin{description}
\item[$\boldsymbol{k=0}$] In this case there is no $SL_2$-fixed point
in $\tilde q^{-1}(\Delta)$, and consequently none on $D_0$. Thus, $X$
must be a cone or $\mathbb Q_3$; see the proof of
corollary~\ref{S2_b21:coprimeness} for this. A contradiction to the
assumption.

\item[$\boldsymbol{k=1}$] In this case $(x,y,z)$ are coordinates for
$\tilde X$. The quotient is terminal iff the weights are of the form
$(a,-a,1)$ (see \cite[sect.~5.1]{Reid87}). Thus, $n'=2$.

\item[$\boldsymbol{k>1}$] Note that there is no $\ZZ_{n'}$-fixed
subspace in $\CZ^3\times\Delta$. In this situation
\cite[Thm. 12]{Mori85} shows that $n'$ must be 4.
\end{description}

Now apply the same argumentation to $D_\infty$ and realize that the
coprimeness assertion of corollary~\ref{S2_b21:coprimeness} is
necessarily violated. This yields the claim.
\end{proof}

\begin{proof}[Proof of~\ref{S2_b21:MainThm} in case (3) of
corollary~\ref{S2_b21:coprimeness}]
As before, set $X^0:=X\setminus D_\infty$ and consider the divisor
$L:=K_{X^0}-D_0$. By adjunction formula, $L|_{D_0}=K_{D_0}$ which has
index $\frac{n}{2}$. Thus, the index of $L$ in $X^0$ is in
$\frac{n}{2}\NZ$. Now perform the cyclic cover associated to $L$
(see~\cite[cor.~1.9]{Reid79} for details): $\gamma : \tilde X^0\rightarrow
X^0$. \name{Stein} factorization gives a diagram 
$$
\begin{CD}
\tilde X^0 @>{\gamma}>> X^0 \\
@V{\tilde q}VV   @VV{q}V \\
\tilde Y @>>> Y
\end{CD}
$$ 
where we can choose $\tilde Y$ to be normal, hence smooth. We are
interested in the preimage of $D_{0}$. First, note that every vector
field on $X^{0}\setminus Sing(X^0)$ can be lifted to $\tilde
X\setminus q^{-1}(Sing(X^{0}))$. Since $\tilde X$ is normal, we obtain
an action of the associated 1-parameter group on $\tilde X$. In
particular, since $S$ is simply connected, $S$ acts on $\tilde X^0$ in a
way that $\gamma$ is equivariant.

As a next step we need to show that $\tilde q^{-1}(0)$ is reduced. By
corollary~\ref{S2_b21:multiplicities}, $q^{-1}(0)$ has multiplicity
$n'=\frac{n}{4}$. On the other hand, generic $(\gamma\circ q)$-fibers
have at least $\frac{n}{4}$ components. This is due to the fact that
$D_\alpha\cong \Ptwo$, where $SL_2$ acts without fixed point, admits
only $\Sigma_0$ as a connected $S$-equivariant cover. Consequence:
$\tilde{q}^{-1}(0)=\gamma^{-1}(D_0)$ is reduced and isomorphic to
$\mathbb F^a_2$. Now apply the argumentation from the proof in case
(2).
\end{proof}

\end{document}